\documentclass{article}

\usepackage{amsmath,amssymb}

\newcommand\zero{{\bf 0}}
\newcommand\bfc{{\bf c}}
\newcommand\bfv{{\bf v}}
\newcommand\x{{\bf x}}
\newcommand{\bfla}{{\boldsymbol{\lambda}}}
\newcommand{\bfmu}{{\boldsymbol{\mu}}}
\newcommand{\cc}{{\mathbb C}}
\newcommand{\qq}{{\mathbb Q}}

\newtheorem{theorem}{Theorem}[section]
\newtheorem{lemma}[theorem]{Lemma}
\newtheorem{proposition}[theorem]{Proposition}

\newcommand\qed{{\hspace*{\fill}$\Box$\vskip12pt plus 1pt}}

\begin{document}

\title{Sweeping Algebraic Curves for Singular Solutions\thanks{This 
material is based upon work supported by the National Science Foundation
under Grant No.\ 0713018.}}

\author{Kathy Piret\thanks{Department of Mathematics, Statistics,
and Computer Science,
University of Illinois at Chicago,
851 South Morgan (M/C 249),
Chicago, IL 60607-7045, USA.
{\tt email:} kathywang@math.uic.edu}
 \and Jan Verschelde\thanks{Department of Mathematics, Statistics,
and Computer Science,
University of Illinois at Chicago,
851 South Morgan (M/C 249),
Chicago, IL 60607-7045, USA.
{\tt email:} jan@math.uic.edu
{\tt URL:} http://www.math.uic.edu/{\~{}}jan}
}

\date{30 September 2008}

\maketitle

\begin{abstract}
Many problems give rise to polynomial systems.  These systems often
have several parameters and we are interested to study how the 
solutions vary when we change the values for the parameters.
Using predictor-corrector methods we track the solution paths.
A point along a solution path is critical when the Jacobian matrix
is rank deficient.
The simplest case of quadratic turning points is well understood,
but these methods no longer work for general types of singularities.
In order not to miss any singular solutions along a path
we propose to monitor the determinant of the Jacobian matrix.
We examine the operation range of deflation and relate the
effectiveness of deflation to the winding number.
Computational experiments on  systems coming from different
application fields are presented.



\end{abstract}

\section{Introduction}

We consider systems $f(\x,\bfla) = \zero$ 
of $N$ polynomial equations $f = (f_1,f_2,\ldots,f_N)$ 
in $n$ variables $\x = (x_1,x_2,\ldots,x_n)$
and $m$ parameters $\bfla = (\lambda_1,\lambda_2,\ldots,\lambda_m)$,
with complex coefficients: $f_k \in \cc[\x,\bfla]$, $k=1,2,\ldots,N$.
In this paper we restrict to isolated solutions
so we assume $N \geq n$ and most often $N = n$. 
For random choices of the parameters~$\bfla$, we expect all solutions 
to~$f(\x,\bfla) = \zero$ to be isolated and well conditioned.  
Even though the total number of solutions may grow exponentially
in the dimensions and degrees of the system, numerical homotopy
continuation algorithms are efficient to enumerate all solutions.

In many applications one wants to know for which values
of the parameters~$\bfla$, the corresponding values of~$\x$
satisfying $f(\x,\bfla) = \zero$ give rise to {\em singular} solutions.
For values of~$\bfla$, the solution~$\x$ is singular if the Jacobian
matrix $A(\x)$, where $A_{ij} = \frac{\partial f_i}{\partial x_j}$,
for $i=1,2,\ldots,N$ and $j=1,2,\ldots,n$, has rank strictly less
than~$n$ when evaluated at~$\x$.  A value for $\bfla$ is {\em critical} 
if some corresponding values for~$\x$ are singular solutions.
To locate all critical values, we may solve an augmented system:
\begin{equation} \label{eqjacocrit}
    F(\x,\bfla,\bfmu) = 
    \left\{
       \begin{array}{lcr}
          f(\x,\bfla) & = & \zero \\
          A(\x,\bfla) \bfmu & = & \zero \\
          \multicolumn{1}{r}{\bfc^T \bfmu} & = & 1
       \end{array}
    \right.
\end{equation}
where $\bfmu$ consists of $n$ additional multiplier variables
and $\bfc \in \cc^n$ is a tuple of $n$ random complex numbers.
The condition $\bfc^T \bfmu = 1$ implies that we look for 
solutions $(\x,\bfla,\bfmu)$ to express $\zero$ as a nonzero
linear combination of the columns of the Jacobian matrix~$A(\x,\bfla)$.
The system in~(\ref{eqjacocrit}) is an example of
{\em the Jacobian criterion} to find critical values.

Solving $F(\x,\bfla,\bfmu) = \zero$ is an effective global method
to locate all critical values~$\bfla$, but because the input size 
of $F$ is more than double the size of the original system~$f$,
this global method is often too expensive --- its underlying
complexity is that of the discriminant variety.
In this paper we focus on the local approach:
given a sufficiently generic value for~$\bfla$,
we start at the corresponding solutions for $\x$ and trace the
algebraic curves as we {\em sweep} through the parameter space.

For introductions to homotopy continuation methods
specific for polynomial systems, we recommend~\cite{Li03}, 
\cite{Mor87}, and~\cite{SW05}.
The books~\cite{AG03}, \cite{Gov00}, and~\cite{Mei00}
provide introductions to path following methods applied to
general nonlinear systems and systems of differential equations.
The authors of~\cite{DLA06} study polynomial differential systems
in the real plane and developed software to draw phase portraits.
Computer algebra is used to compute all singularities, but it is
noted in~\cite{DLA06} that for high degrees this can take a long time.
Recent related symbolic methods are described in~\cite{LR07}
and~\cite{SS04}.

In most applications, one is mainly interested in real solutions.
However, a complex solution curve of a polynomial system may have
isolated real solutions.  Such a real solution on a complex curve
will be isolated in the real space and will manifest itself as
a singular solution on the curve.  The methods in~\cite{LBSW07}
rely on a global application of a deflation operator to locate
real solutions on a curve in complex space while the sweeping homotopies
defined in this paper offer a local approach for this problem.

The contributions of this paper are twofold.  To detect general
types of singularities along a solution path we first describe an
algorithm to monitor the determinant of the Jacobian matrix.
Secondly, we investigate the effectiveness of deflation to
accurately locate the detected isolated singular solution.
These two contributions are outlined in sections 4 and~5.
In sections~2 and~3 we first define sweeping homotopies and
describe the problem statements.  This paper ends with 
a report on our computational experiments.

\section{Sweeping Homotopies}

The system $f(\x,\bfla) = \zero$ defines already a homotopy.
We call it a {\em natural parameter homotopy} because the
parameters $\bfla$ appear naturally.  To track the solution
paths $\x(\bfla)$ with pseudo arc length continuation we first compute
a tangent vector $\bfv = (\bfv_\x , \bfv_\bfla)^T$
at the current point $(\x_0,\bfla_0) \in \cc^n \times \cc^m$:
\begin{equation}
   \left(
      \begin{array}{cc}
         \frac{\partial f}{\partial \x}(\x_0,\bfla_0)
         & \frac{\partial f}{\partial \bfla}(\x_0,\bfla_0)
      \end{array}
   \right)
   \left(
      \begin{array}{c}
          \bfv_\x \\ \bfv_\bfla
      \end{array}
   \right), \quad || \bfv || = 1.
\end{equation}
At a regular solution on an algebraic curve defined by
a complete intersection where $N = n$ and $m = 1$, 
the unit tangent vector is defined uniquely up to orientation.
After selecting an orientation, a prediction for the next point
is then $(\x_1,\bfla_1) = (\x_0,\bfla_0) + h (\bfv_\x, \bfv_\bfla)$
where $h$ is the step size.  Using interval methods as in~\cite{KX94}
to control the step size, one can rigorously prevent jumping from one
path to another.  After the prediction step, Newton's method is applied
to correct the predicted solution back to the solution curve.
Algorithms for the adaptive use of multiprecision 
arithmetic during path following are proposed in~\cite{BHSW08}
and~\cite{JZ05}.

In our sweep we use an {\em artificial parameter homotopy},
introducing a new artificial parameter~$t$ to make a convex combination
between two given sets of parameter values~$\bfla_0$ and $\bfla_1$.
Given $f(\x,\bfla) = \zero$, $\bfla_0$ and $\bfla_1$, 
a {\em sweeping homotopy} is defined as
\begin{equation}
   h(\x,\bfla,t) = 
   \left\{
       \begin{array}{r}
          f(\x,\bfla) = \zero \\
          (1-t) (\bfla - \bfla_0) + t (\bfla - \bfla_1) = \zero
       \end{array}
   \right.
   \quad t \in [0,1].
\end{equation}
For $t = 0$, we start a solutions of $f(\x,\bfla_0) = \zero$ and
as $t$ moves to 1, we sweep to the solutions of $f(\x,\bfla_1) = \zero$.
We use the same type of pseudo arc length continuation as above,
with $t$ as an added parameter, enforcing the orientation of the
tangent vector so that $t$ always strictly increases in value.

The distinction between natural and artificial parameter homotopies
has profound consequences for the treatment of singularities.
Consider for example $f(x,\lambda) = x^2 + \lambda^2 - 1 = 0$.
Viewing $f(x,\lambda) = 0$ naturally, we recognize the equation
of the (real) unit circle in the plane.  Tracking the curve
$x(\lambda)$ as defined by the natural homotopy $f(x,\lambda) = 0$
is then simply tracing the circle, either clockwise or counterclockwise.
Consider picking $\lambda_0 = 0$ and $\lambda_1 = 2$ as start
and target values in a sweeping homotopy:
\begin{equation}
   h(x,\lambda,t) = 
   \left\{
       \begin{array}{r}
          x^2 + \lambda^2 - 1 = 0 \\
          (1-t) \lambda + t (\lambda - 2) = 0
       \end{array}
   \right.
   \quad t \in [0,1].
\end{equation}
Forcing the orientation of the tangent vector to the path
$(x(t),\lambda(t),t)$ so $t$ is strictly increasing leads
to a quadratic turning point for $t = 0.5$.  At that point,
the two real paths turn into the complex plane.
While real homotopies (i.e.: homotopies with all coefficients real)
for solving polynomial systems lead to singular points along the
solution paths, as shown in~\cite{LW93}, generically, only a finite
number of quadratic turning points occur.
However, for our problem, perhaps we might assume a generic choice 
of the values for the parameters~$\bfla_0$ at the start of the sweep, 
but even that is insufficient to exclude general types of singularities.

\section{Detection and Location of Singularities}

Given a sweeping homotopy $h(\x,\bfla,t) = \zero$
and a start solution $(\x_0,\bfla_0,0)$, 
our problem is then to detect and locate all singular
solutions along the path as $t$ advances from 0 to~1.

For the most common type of singularity, the quadratic turning
points, the detection and location of singularities along the
solution paths is done as follows:
\begin{description}
\item[Detection: via the orientation of the tangent vectors]  
     The tangent vector $\bfv$ has three components
     $\bfv = (\bfv_\x, \bfv_\bfla, \bfv_t)$ and each time
     we force its orientation so $\bfv_t > 0$.
     If $t_1$ is before and $t_1$ after a turning point,
     forcing $\bfv_t > 0$ at $t_2$ will lead to a change in 
     the angle between the corresponding tangent vectors
     $\bfv(t_1)$ and $\bfv(t_2)$ so that its 
     inner product $\langle \bfv(t_1) , \bfv(t_2) \rangle < 0$.
\item[Location: shooting method for the step size]
     Once we have two consecutive tangent vectors
     $\bfv(t_1)$ and $\bfv(t_2)$ for which 
     $\langle \bfv(t_1) , \bfv(t_2) \rangle < 0$ we look to find $h$
     so that $\langle \bfv(t_1) , \bfv(t_1 + h) \rangle = 0$.
     To find $h$ we apply a shooting method.  We found the
     description in~\cite{LZC92} very clear and useful.
\end{description}
The two tasks {\em detection} and {\em location} of a singularity
along a solution path defined by a sweeping homotopy are more
accurately described in the input/output statements 
in Table~\ref{tabdetprob} and Table~\ref{tablocprob}.

\begin{table}[h]
\begin{center}
\begin{tabular}{|rl|} \hline
 \multicolumn{2}{|c|}{Input/Output specification for Detection Problem}
\\ \hline
  Input : & $h(\x,\bfla,t) = \zero$ sweeping homotopy, \\
          & $(\x_0,\bfla_0,0)$ a start solution. \\
 Output : & solutions 
           $(\x(t_1),\bfla(t_1),t_1)$ and $(\x(t_2),\bfla(t_2),t_2)$, \\
          & where the interval $(t_1,t_2)$ contains a singularity. \\ \hline
\end{tabular}
\caption{The Detection Problem: detect singularities along a path.}
\label{tabdetprob}
\end{center}
\end{table}

The output of the Detection Problem will be empty if there are no
singularities for all $t \in [0,1]$.
The main difficulty is to distinguish this case
from the case where the solutions paths run straight through a
multiple solution.  If the paths are straight, the path tracker will
not slow down and overshoot the singular solution.
The output of the Detection Problem will be incomplete in the 
case where paths near a singular solution are very hard to follow.
An example of that case is when a severe drop in the rank of the 
Jacobian matrix causes Newton's method to fail.  
The output in such a case then only consist of $t_1$
and its corresponding solution, even as we have a good guess
for $t_2$, the solution corresponding to $t_2$ will be missing.

The difficulties in both problems are often complementary.
On the one hand, singularities that are easy to detect 
(because Newton's method fails), are usually hard to 
locate (again for the same reason).
On the other hand, multiple solutions for which Newton's method
converges are hard to detect along a path.

The output of the Detection Problem determines the input for the
Location Problem, specified in Table~\ref{tablocprob}.

\begin{table}[h]
\begin{center}
\begin{tabular}{|rl|} \hline
 \multicolumn{2}{|c|}{Input/Output specification for Location Problem}
\\ \hline
  Input : & $h(\x,\bfla,t) = \zero$ sweeping homotopy, \\
          & $(\x(t_1),\bfla(t_1),t_1)$ and $(\x(t_2),\bfla(t_2),t_2)$
            are solutions \\
          & where the interval $(t_1,t_2)$ contains a singularity. \\
 Output : & $(\x^*,\bfla^*,t^*)$ accurate singular solution
            of $h(\x,\bfla,t) = \zero$. \\ \hline
\end{tabular}
\caption{The Location Problem: locate accurately a singularity along a path.}
\label{tablocprob}
\end{center}
\end{table}

The quality of the input to the Location Problem will determine
the difficulty of the Detection Problem.  If the corresponding
solutions at $t_1$ and $t_2$ are accurate and not too far apart from
each other, then solving the Location Problem will be much easier
than when the corresponding solution for $t_1$ is inaccurate and
the solution corresponding to $t_2$ is missing.

In the next two sections we outline our algorithms for the
detection and location problem.

\section{Detecting Singularities along a Path}

The main difficulty in the detection problem is that the path tracker
may not slow down when approaching a well conditioned multiple solution.
Therefore, we first look for a criterion to decrease the step size
and to back up towards previously computed values.

The orientation of the tangent offers a clear criterion 
for quadratic turning points, but is no longer useful
when sweeping paths that do not turn at a singularity.
Monitoring the signs of the eigenvalues of the Jacobian matrix
also captures many types of singular solutions, see e.g.~\cite{GMS96},
but we have encountered cases -- see the applications section below -- 
where the eigenvalues do not change signs when passing through
the singular solution and where thus the determinant of the
Jacobian matrix stays monotone of the same sign, only touching
zero at the singular solutions.  Our experiences have led us to opt 
for the determinant of the Jacobian matrix as the main criterion.
This determinant is obtained as a relatively easy byproduct of 
the application of Newton's method. 

We can compute the determinant only we have as many equations~$N$
as unknowns~$n$.  In case~$N > n$ we can locally make the system
square either by adding random combinations of the extra $N-n$
polynomials to first $n$ polynomials, or by adding $N-n$ slack
variables to all polynomials, see~\cite[\S 13.5]{SW05}.

To detect singularities, we keep a window of three consecutive
values for the artificial parameter $t$: $t_1 < t_2 < t_3$, 
along with the values of the determinants $d_1$, $d_2$, and $d_3$
of the Jacobian matrix at the corresponding solutions
$\x(t_1)$, $\x(t_2)$, and $\x(t_3)$. 
If there are any sign changes in the determinants, 
then the detection problem is solved.
Otherwise, we compute an interpolating parabola $p(t)$
so that $p(t_k) = d_k$, for $k=1,2,3$.
If all determinants are positive, we compute the minimum of~$p$.
If all determinants are negative, we compute the maximum of~$p$.
The distinction on the sign of the determinant only matters for
{\em real} solution paths.  For solution paths in complex space,
we monitor the modulus of the determinant, 
i.e.: $d_k = |\det(A(\x(t_k),\bfla(t_k),t_k)|$, $k=1,2,3$.
An explicit criterion is given in~(\ref{eqcriterion}) below.

\begin{lemma} \label{lemcriterion}
Let $d_1$, $d_2$, and $d_3$ correspond to
the three consecutive values for $t$: $t_1 < t_2 < t_3$.  Then
\begin{equation} \label{eqcriterion}
   z = \frac{~~t_1^2 (d_3 - d_2) + t_2^2 (d_1 - d_3) + t_3^2 (d_2 - d_1)}
     {2 \left( d_1   (t_2 - t_3) + d_2   (t_3 - t_1) + d_3   (t_1 - t_2)
       \right)}
\end{equation}
is a critical value for the interpolating parabola
through the points $(t_1,d_1)$, $(t_2,d_2)$ and~$(t_3,d_3)$.
\end{lemma}
{\em Proof.}  The Lagrange form of the parabola interpolating 
through the points $(t_1,d_1)$, $(t_2,d_2)$ and~$(t_3,d_3)$ is
\begin{equation}
   p(t) = \frac{(t-t_2)(t-t_3)}{(t_1-t_2)(t_1-t_3)} d_1
        + \frac{(t-t_1)(t-t_3)}{(t_2-t_1)(t_2-t_3)} d_2
        + \frac{(t-t_1)(t-t_2)}{(t_3-t_1)(t_3-t_2)} d_3.
\end{equation}
Computing the derivative $p'(t) = 0$ and solving for $z$
in $p'(z) = 0$ yields~(\ref{eqcriterion}).~\qed

\begin{lemma} \label{lemcostcrit}
The cost of evaluating~{\rm (\ref{eqcriterion})} along
a path in $\cc^n$ is $O(n)$.
\end{lemma}
{\em Proof.}  Tracking a path, using Newton's method
as a corrector, the Jacobian matrix is evaluated and decomposed 
when solving a linear system to obtain the next iteration of
Newton's method.  Given an LU decomposition of the Jacobian matrix,
the extra cost of computing the determinant is just a product of
$n$ numbers along the diagonal of one of the factors $L$ or~$U$.
Once the points $(t_1,d_1)$, $(t_2,d_2)$ and~$(t_3,d_3)$ are given,
evaluating~(\ref{eqcriterion}) requires a constant number of
arithmetical operations.  Thus the cost of evaluating~(\ref{eqcriterion})
along a path in $\cc^n$ is $O(n)$.~\qed

Lemmas~\ref{lemcriterion} and~\ref{lemcostcrit} ensure that we have
an explicit and efficient formula to evaluate along a path.
Taking only the linear algebra operations into account,
one iteration of Newton step costs at least~$O(n^3)$ which
already dominates the~$O(n)$ cost of evaluating the formula.
To show the effectiveness of~(\ref{eqcriterion}),
we exploit the algebraic properties of our problem.
In particular, we use Puiseux expansions~(\cite{dJP00}, \cite{Wal50}).

\begin{lemma} \label{lempuiseux}
Assume the sweeping homotopy $h(\x,\bfla,t) = \zero$ has
an isolated finite singularity for $t = t_*$.
Then the determinant of the Jacobian matrix of $A(\x(t),\bfla(t))$
for $t$ sufficiently close to~$t_*$ equals $(t-t_*)^p$, 
for some fractional power~$p$.
\end{lemma}
{\em Proof.}  A Puiseux expansion at an isolated singular solution of a 
polynomial system is a fractional power series,
i.e.: with rational numbers for the powers in the series.
In particular, for $t \approx t_*$ we may write
$x_k(t) = \alpha_k (t-t_*)^{a_k}(1 + O(t))$, $k=1,2,\ldots,n$,
where $\alpha_k \in \cc \setminus \{ 0 \}$ and $a_k \in \qq$.
Because we assumed the singularity to be finite
(i.e.: not at infinity), we have that $a_k > 0$ for all~$k$.
By definition of the sweeping homotopy,
$\bfla(t)$ is linear in $t$ and thus also in $t - t_*$.

Substituting $x_k(t) = c_k (t-t_*)^{a_k}(1 + O(t))$
and the linear expression for $\lambda_k(t)$
into the Jacobian matrix $A(\x(t),\bfla(t),t)$, we view $A$ as $A(t)$,
as a matrix of polynomials with fractional powers in~$t-t_*$.
So also its determinant, $\det(A(t))$ is a polynomial in~$t-t_*$.
Ignoring higher order terms, we let $p$ be the smallest 
power of $t-t_*$ in~$\det(A(t))$.
For $t \approx t_*$, we then have
$\det(A(t)) \approx (t-t_*)^p$.~\qed

If the power~$p$ of Lemma~\ref{lempuiseux} is a natural number,
then the interpolating parabola of Lemma~\ref{lemcriterion}
will locally resemble very well the determinant itself
and its critical value will be close to~$t_*$.
For fractional powers of~$p$, the determinant does not look
like a polynomial.  For example, consider $t_* = 0$ and $p = 1/2$,
then for $t > 0$, $\det(A(t)) = \sqrt{t}$
and for $t < 0$, $\det(A(t)) = \sqrt{-t}$.

\begin{theorem} \label{thecriterion}
Assume the sweeping homotopy $h(\x,\bfla,t) = \zero$
has exactly one isolated finite singular solution 
at $t_* \in [t_1,t_3]$.
Then for any~$t_2 \in (t_1,t_3)$ and~$z$ as 
in~{\rm (\ref{eqcriterion})}: $z \in [t_1,t_3]$.
\end{theorem}
{\em Proof.}  Without loss of generality, consider~$t_* = 0$.
If~$p$ of Lemma~\ref{lempuiseux} is a natural number,
then for $t \approx 0$: $\det(A(t)) \approx t^p$
and the interpolating parabola of Lemma~\ref{eqcriterion}
will have its critical value~$z$ in $[t_1,t_3]$.
Even if~$p$ is not a natural number and a fraction,
for all $t_2 \in (t_1,t_3)$ the value of $\det(A(t_2))$
will be smaller than $\det(A(t_1))$ and $\det(A(t_3))$ 
so the interpolating parabola of Lemma~\ref{eqcriterion}
will have the right concavity and thus its critical value 
$z$ will also be in~$[t_1,t_3]$.~\qed

Once the formula~(\ref{eqcriterion}) yields 
a value~$z \in [t_1,t_3]$ we then have a candidate singularity
at~$t^*$ and we need to locate it accurately.
The case where the determinant of the Jacobian matrix has a local
optimal value in~$[t_1,t_3]$ but no root is captured by the application of
a minimization algorithm on the determinant, viewed as a function in~$t$.
Although we could further apply parabolic interpolation and
use~$z$ of~(\ref{eqcriterion}) as the next value for~$t$,
for fractional powers, we recommend the golden section search 
method~\cite{NW99}.
The golden section search method requires an optimal number of
function evaluation and is guaranteed to find the optimum if
the function is unimodal over the interval.

The key assumption of Theorem~\ref{thecriterion} is
that there is only one singularity in the interval.
For a reliable implementation of this algorithm,
we must relate the step size $h$ of the path tracker
to the distance $\delta$ between two singular solutions.
If $h$ is sufficiently smaller than $\delta$,
then it is safe to assume that the determinant of the Jacobian matrix
will behave as a unimodal function between three consecutive 
predictor-corrector steps.
For polynomial systems, the total degree $D$ is the product of the
degrees of the polynomials in the system.  Following B\'ezout's 
theorem it is a crude upper bound on the number of isolated solutions
and therefore also on the total number of singularities.
Assuming a uniform distribution of the singularities,
a pessimistic lower bound on the step size $h$ is $1/D^2$.

\section{Locating Singularities along a Path}

If a singular solution at $t^*$ is hard to detect,
then for almost all $t$ close to $t^*$
the Jacobian matrix is sufficiently well conditioned for
Newton's method to converge well.  In that sense, getting
close enough to singularity to locate it with sufficient
accuracy is then no problem.
Thus then the main difficulty with the location problem
occurs when Newton's method fails.

The solution to the detection problem
has made the location problem similar 
to an {\em endgame}~\cite{MSW92b} (see also~\cite{SWS96}).
In this section we discuss the effectiveness of deflation
--- the idea to apply deflation to locate singular solutions
of polynomial systems occurred first in~\cite{OWM83},
see~\cite{Lec02} for a symbolic deflation method ---
to locate general types of isolated singular solutions
in the context of a sweeping homotopy.

The deflation operator works {\em locally} starting at
an approximation for a singular solution for which the
Jacobian matrix has numerical rank is~$R$.
Numerical rank revealing algorithms can be found in~\cite{LZ05}.
Then $R+1$ multiplier variables~$\bfmu$ are used.
To reduce to the corank one case, we multiply the 
Jacobian matrix with a random matrix $B$.
This matrix~$B$ has  $R+1$ columns and as many rows 
as the columns of the Jacobian matrix.
Then we apply Newton's method on the system

\begin{equation} \label{eqdeflation}
    E(\x,\bfla,\bfmu) = 
    \left\{
       \begin{array}{lcr}
          f(\x,\bfla) & = & \zero \\
          A(\x,\bfla) B \bfmu & = & \zero \\
          \multicolumn{1}{r}{\bfc^T \bfmu} & = & 1
       \end{array}
    \right.
\end{equation}
The system $E(\x,\bfla,\bfmu) = \zero$ looks very similar
to the augmented system of the Jacobian criterion~(\ref{eqjacocrit}),
with the addition of the numerical rank as extra local information.
One application of the deflation operator may not be enough to
completely recondition the isolated singularity and we have to
apply deflation recursively.
As proven in~\cite{DZ05} and~\cite{LVZ06}, the number of deflations 
needed to restore the quadratic convergence of Newton's method is strictly
less than the multiplicity of the singular solution.

The term {\em endgame operation range} was coined in~\cite{MSW92b}.
In general, this endgame operation range is the range for which the
endgame techniques are effective.  If we use for example extrapolation
methods, then we need on the one hand take samples along the path close
enough to the singularity.  On the other hand, if we get too close to
the singularity, the iterations produced by Newton's method will be
too inaccurate for the extrapolation.  If we can adjust the working
precision of our calculations, then we can guarantee that the endgame
operation range is nonempty.

The idea of deflation is to consider in addition of the original
polynomials also the derivatives.  If on the one hand, we are too far from 
the singularity, then adding the derivatives at the current approximation 
may lead to an inconsistent problem and lead to divergence.  On the other
hand, getting close enough to the singular solution may no longer be
possible by the plain application of Newton's method.

For deflation, one critical factor in its endgame operation range 
is {\em the winding number}.  The winding number occurs as the denominator
in the fractional power (or Puiseux) expansion of the solution path at
the singular solution.  The multiplicity of the solution bounds the
winding number.  The higher the winding number, the harder it could be
for the derivatives to vanish in the proximity of the singularity.
Proposition~\ref{proprangedeflation} formalizes
the relationship between the winding number and 
the endgame operation range for deflation.

\begin{proposition} \label{proprangedeflation}
Let $h(\x,\bfla,t) = \zero$ be a sweeping homotopy with an isolated
finite singular solution for $t = t_*$ with winding number~$\omega$.
If for some component $k$:
$h_k(\x(t),\bfla(t),t)$ is $O(t)$ for $t \approx t_*$,
then $\frac{\partial h_k}{\partial x_j} h(\x(t),\bfla(t),t)$ could 
in the worst case be~$O(t^{1/\omega})$.
\end{proposition}
{\em Proof.}  Without loss of generality we may assume that $t_* = 0$,
$\bfla(t_*) = \zero$ and $\x(t_*) = \zero$.  Expanding the solutions
$\x(t)$ at $t_* = 0$ leads to fractional power series
$x_i(t) = c_i t^{v_i/\omega}(1 + O(t))$, for $i=1,2,\ldots,n$,
where $c_i \in \cc \setminus \{ 0 \}$ and $v_i$ is a natural number.
By definition of the sweeping homotopy, the relation between $\bfla$
and~$t$ is simply linear and it is straightforward to express~$\bfla(t)$
as a linear function of~$t$.

Substituting $\bfla(t)$ and the power series for $\x(t)$ in $h_k$
we obtain $h_k(\x(t),\bfla(t),t) = \gamma t^p(1+O(t))$ for some nonzero
complex constant $\gamma$ and some power~$p$.
Since we assumed a finite singular solution: $p \geq 1$.
Similarly, for a derivative $\frac{\partial h_k}{\partial x_j}$
we obtain $\frac{\partial h_k}{\partial x_j}(\x(t),\bfla(t),t) 
= \delta t^q (1 + O(t))$ for some nonzero complex constant~$\delta$
and some power~$q$.
Because not all derivatives will vanish at the singular solution,
suppose $k$ and $i$ are such that $q$ is the lowest positive exponent.
Take then $p = (\omega+1)/\omega$ and let $q = p-1$.~\qed

As a practical result of Proposition~\ref{proprangedeflation}
we may make some pessimistic predictions on the endgame operation
range of deflation.  For example, if $\omega = 4$ and we need
the residual of the derivatives to be about $10^{-2}$, then
the residual of the approximation must be about $10^{-8}$.

Also the numerator of the exponent in the leading term of the
Puiseux series plays an important role as it determines how sharp
the curve bends as $t$ approaches~$t_*$.  Even with a high winding
number, extrapolation methods will be effective for low numerators,
but as the numerator is close to $\omega$ itself, then the solution
curve will appear to be linear unless we get really close to~$t_*$.

To estimate the winding number, Richardson extrapolation
(see e.g.:~\cite{BZ91}) cannot be applied directly because
the exponents in the power series are unknown --- we refer
to~\cite{CP89} for general extrapolation methods for unknown exponents.
In~\cite{HV98}, extrapolation methods to estimate the winding
number for diverging solution paths were developed.
An algorithm to predict the order of the deflation was
recently presented in~\cite{LVZ08}.

\section{Computational Experiments}

We have implemented the detection criterion 
in the publicly available open source software PHCpack~\cite{Ver99}
and applied it to three polynomial systems,
coming from different application fields
and documented in the literature~\cite{EM99}, \cite{Noo89}, \cite{WW05}.

\subsection{a system from molecular configurations}

The following system occurs in~\cite{EM99}:
\begin{equation}
   f(\x,\lambda) =
   \left\{
     \begin{array}{l}
        \frac{1}{2}(x_2^2 + 4x_2x_3 + x_3^2) + \lambda (x_2^2 x_3^2 - 1) = 0 \\
   \vspace{-2mm} \\
        \frac{1}{2}(x_3^2 + 4x_3x_1 + x_1^2) + \lambda (x_3^2 x_1^2 - 1) = 0 \\
   \vspace{-2mm} \\
        \frac{1}{2}(x_1^2 + 4x_1x_2 + x_2^2) + \lambda (x_1^2 x_2^2 - 1) = 0.
     \end{array}
  \right.
\end{equation}
The system is listed as a nontrivial example in~\cite[pages 391-392]{Ste04}.
However, the system is small enough for global methods.
The Jacobian criterion~(\ref{eqjacocrit}) gives a system we solved with the 
blackbox solver of PHCpack~\cite{Ver99}.  The system has 54 generic
solutions which can be divided into five groups with the same
$x_1$, $x_2$, $x_3$ and $\lambda$ values. 
The first four groups have the same absolute values 
of $x_1$, $x_2$ and $x_3$ with the natural parameter
$\lambda$ being either $+1.5i$ or $-1.5i$, $i = \sqrt{-1}$.
There are exactly twelve solutions in each of the first four groups. 
The last group corresponds to the approximate value 
$\pm$ 0.866025403780023 as the natural parameter~$\lambda$.
For these two values there are curves of degree six.
In this example, all critical values for $\lambda$ were found
via the Jacobian criterion.

As $\lambda$ approaches zero, the system becomes singular.
At the origin, there is one solution of multiplicity~8 for the
system when the deflation method in PHCpack is applied.
To test our detection algorithm, we consider sweeping $\lambda$
through zero.  The sweeping homotopy is
\begin{equation}
   f(\x,\lambda) =
   \left\{
     \begin{array}{l}
        \frac{1}{2}(x_2^2 + 4x_2x_3 + x_3^2) + \lambda (x_2^2 x_3^2 - 1) = 0 \\
   \vspace{-2mm} \\
        \frac{1}{2}(x_3^2 + 4x_3x_1 + x_1^2) + \lambda (x_3^2 x_1^2 - 1) = 0 \\
   \vspace{-2mm} \\
        \frac{1}{2}(x_1^2 + 4x_1x_2 + x_2^2) + \lambda (x_1^2 x_2^2 - 1) = 0 \\
   \vspace{-2mm} \\
        (\lambda-1)(1-t) + (\lambda + 1)t = 0.
     \end{array}
  \right.
\end{equation}

As the artificial parameter t goes from 0 to 1, the natural
parameter $\lambda$ is swept from $+1$ to $-1$. 
According to the multihomogenous B\'ezout bound~\cite{SW05},
the permanent of the degree matrix of the system is 16 for nonzero 
values of~$\lambda$.  This bound is sharp, so all solutions in 
a multihomogeneous homotopy converge to finite solutions.  
Among the 16 solutions, four are symmetrical complex conjugate 
solution pairs and four are symmetrical real solution pairs.  
By the symmetry, the solutions break up into orbits of type
$x_1=x_2$, $x_2=x_3$, $x_1=x_3$ and $x_1=x_2=x_3$.
As $\lambda$ is swept from $+1$ to $-1$, starting with start solutions 
at $t = 0$, four real solution paths converge around the origin and
the four complex solution paths diverge. 
If we would like to track the converging complex solutions paths, 
we could set the homotopy to
$(\lambda+1)(1-t) + (\lambda - 1)t = 0$ such that the four complex
solution paths converge around the origin and the four real
solution paths diverge. The special value zero for the natural
parameter $\lambda$ is found by the sweep as the tangent flips.
A solution of multiplicity~8 is found at the origin.

Our new detection algorithm is needed to detect the singularities
at $\lambda$ = $\pm$0.866025403780023 for which there are
curves of degree~6.  Because even close to this
critical value, the solutions are still relatively well conditioned,
monitoring just the orientation of the tangent is insufficient.

\subsection{modeling neural networks}

Families of polynomial systems often not only depend on parameters,
but also the dimension may scale.  Our next example originates
from~\cite{Noo89}.  For $n = 3$, an example of a system in this
family is
\begin{equation}
   f(\x,\lambda) =
   \left\{
      \begin{array}{l}
         x_1 x_2^2 + x_1 x_3^2 - \lambda x_1 + 1 = 0 \\
         x_2 x_1^2 + x_2 x_3^2 - \lambda x_2 + 1 = 0 \\
         x_3 x_1^2 + x_3 x_2^2 - \lambda x_3 + 1 = 0. \\
      \end{array}
   \right.
\end{equation}
The application of the Jacobian criterion in~(\ref{eqjacocrit}), 
results in a 7-by-7 system with 54 regular solutions. 
Critical values for~$\lambda$ found among these solutions are
0, and the approximations 1.88988157484231, 3.61703146124952, 
2.38110157795230, and $-0.414704714645147$. 
As the dimension~$n$ grows, the plain application of 
the Jacobian criterion will quickly lead to an intractable problem, 
whereas the complexity of tracking one solution path scales much better.

The singular solutions for the critical values corresponding
to~$\lambda = 0$ were the hardest to detect and stimulated
the development of our detection algorithm.  
The homotopy which sweeps $\lambda$ through zero is
\begin{equation}
   f(\x,\lambda) =
   \left\{
      \begin{array}{r}
         x_1 x_2^2 + x_1 x_3^2 - \lambda x_1 + 1 = 0~ \\
         x_2 x_1^2 + x_2 x_3^2 - \lambda x_2 + 1 = 0~ \\
         x_3 x_1^2 + x_3 x_2^2 - \lambda x_3 + 1 = 0~ \\
         (\lambda + 0.1)(1-t) + (\lambda - 0.1)t = 0. \\
      \end{array}
   \right.
\end{equation}
As the artificial parameter t goes from 0 to 1, the natural
parameter $\lambda$ is swept from $-0.1$ to $0.1$. 
Passing through $\lambda = 0$, the tangent vector does not flip back,
the determinant does not change sign and comparing the signs of
eigenvalues for $\lambda < 0$ and $\lambda > 0$ does not reveal anything.
Without our detection algorithm, the path tracker will not back up and the
solution of multiplicity four for $\lambda = 0$ would
go undetected.  For general values of~$n$, the solution
corresponding to~$\lambda = 0$ has multiplicity~$n+1$.

\subsection{a symmetric Stewart-Gough platform}

A Stewart-Gough platform consists of two plates connected by six legs.
One plate is fixed (the base plate) while the other plate (the top
plate) moves as the leg lenghts change.
These platforms are use used for example in flight simulators.
At a singularity the trajectory of the top plate
is no longer uniquely defined.
In our experiments, we follow~\cite{WW05} where the equations
for a symmetric Stewart-Gough platform are:

\begin{equation}
   f(\x,l_1) =
   \left\{
     \begin{array}{l}
        (x_i-x_{i0})^2+(y_i - y_{i0})^2 + z_i^2 - l_i^2 = 0, i= 1,2,...,6 \\
        (x_2-x_1)^2 + (y_2 - y_1)^2
	(z_2 - z_1)^2 - 2R_1^2(1-\cos(\alpha_1)) = 0 \\
        (x_1-x_0)^2 + (y_1 - y_0)^2 + (z_1 - z_0)^2 - R_1^2 = 0 \\
        (x_2-x_0)^2 + (y_2 - y_0)^2 + (z_2 - z_0)^2 - R_1^2 = 0 \\
     \end{array}
  \right.
\end{equation}
  for
\begin{equation}
   \left\{
     \begin{array}{l}
        x_i = w_1 x_0 + w_2^{m_1}w_3^{m_2}x_1 + w_2^{m_2}w_3^{m_1}x_2 \\
        y_i = w_1 y_0 + w_2^{m_1}w_3^{m_2}y_1 + w_2^{m_2}w_3^{m_1}y_2 \\
        z_i = w_1 z_0 + w_2^{m_1}w_3^{m_2}z_1 + w_2^{m_2}w_3^{m_1}z_2  \\
     \end{array}
  \right.
\end{equation}
 where
 \begin{equation}
   \left\{
     \begin{array}{l}
        w1 = \frac{3\sin(\alpha_1) + (-1)^m \sqrt3(\cos(\alpha_1) - 1)}
                  {2\sin(\alpha_1)} \\
   \vspace{-2mm} \\
        w2 = \frac{-\sin(\alpha_1 
              - (-1)^m \sqrt3\cos(\alpha_1)}{2\sin(\alpha_1)} \\
   \vspace{-2mm} \\
        w3 = \frac{(-1)^m \sqrt3}{2\sin(\alpha_1)}.
     \end{array}
  \right.
\quad
   \left\{
     \begin{array}{l}
        m =0, \ {\rm for} \ i =3,6, \\
        m =1, \ {\rm for} \ i =4,5  \\
        m_1 =0,m_2=1,\  {\rm for} \  i =3,5, \\
        m_1 =0,m_2=1,\  {\rm for} \  i =4,6  \\
     \end{array}
  \right.
\end{equation}

The polynomial system has three fixed parameters: $\alpha_1$, $\alpha_2$, 
and $R_1$ which determine the configuration of the platform.
The angle $\alpha_1$ is the relative angle between two the triangles 
connecting the joints in the moving top platform, 
while $\alpha_2$ is the relative angle between two triangles connecting
the joints in the fixed base platform.
The radius $R_1$ is the radius of joints on the top plate.
As in~\cite{WW05}, we fix the configuration parameters: $R_1 = 1$,
$\alpha_1$ and $\alpha_2$ are respectively 28 and 22 degrees.
Although we could consider the system as depending on six parameters,
the leg lengths $l_i$, $i=1,2,\ldots,6$, for the purpose of simplicity,
we only treat $l_1$ as a natural parameter.

The symmetrical platform gives rise to a system of nine polynomial
equations in nine unknowns $\x = (x_0,y_0,z_0,x_1,y_1,z_1,x_2,y_2,z_2).$
In the application of the Jacobian criterion, we need to solve
a 19-by-19 polynomial system.  Fortunately the system is sparse
and the mixed volume of the tuple of Newton polytopes equals 4,608.
Tracking 4,608 paths yields 256 regular solutions.

Applying our sweep to find critical values is certainly
much less expensive for this system.
Fixing $l_i$ to 1.5, 2.0, and 3.0, we found four special values 
for the natural parameter $l_1$ for each $l_i$ with higher precision
than what was reported in~\cite{WW05}.  In addition, we are able to 
see that $z_{0}$ can be either posive or negative.
When $l_i$ is around 1.003, a multiple singular point occurs at
the origin and $l_i$ approximates the value of $l_1$, the system
becomes a two-parameter problem and requires special care.

\bibliographystyle{plain}

\end{document}